\documentclass{article}
\usepackage{graphicx}

\usepackage{amsmath,amssymb}
\usepackage{amsthm}
\usepackage{amsfonts}
\usepackage{latexsym}
\usepackage{color}
\usepackage{verbatim}
\usepackage{tikz}
\usepackage{xcolor}
\usetikzlibrary{matrix, intersections}
\usepackage[margin=1in]{geometry}

\usepackage{forest}
\usepackage{soul}

\usepackage[utf8]{inputenc}

\theoremstyle{plain}
\newtheorem{theorem}{Theorem}[section] 

\newtheorem*{axiom}{Axiom} 

\theoremstyle{definition}

\theoremstyle{remark}

\title{Is There An Ideal Color Wheel?}
\author{Tejo Madhavarapu, T. Kyle Petersen, and Peter Winkler}

\date{\today}

\begin{document}

\maketitle

\begin{abstract}
The familiar {\em color wheel} is a disk divided into six sectors, colored red, orange, yellow, green, blue, and purple, in
circular order.  Three of the colors can be obtained by blending the colors in the two neighboring sectors.

One might wonder: is there a color wheel in which all six of the sections have this property, without all the sections
being the same color?  We show that the answer is no, not just for the 6-cycle but for any finite connected graph; indeed,
for any finite, strongly connected, edge-weighted digraph.  The result generalizes the ``harmonic lemma" for graphs,
replacing the well-behaved averaging function by paint blending, about which almost nothing is assumed.  Our proof
makes use of a Markov chain stopping rule.
\end{abstract}

\section{A Chromonic Lemma}

\begin{figure}
\centering
\includegraphics[width=6cm]{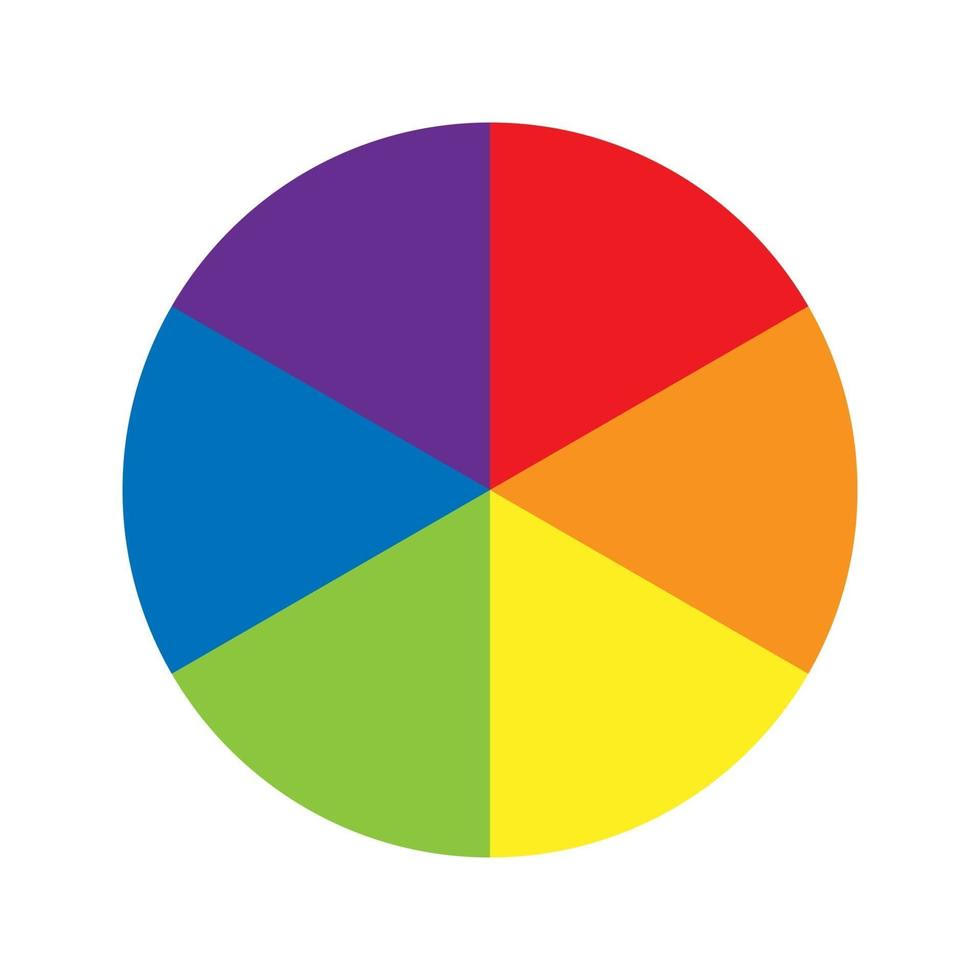}
\caption{The color wheel.}\label{fig:wheel}
\end{figure}

The color wheel shown in Figure \ref{fig:wheel} is taught in elementary school art classes worldwide.
The order in which the colors are portrayed is not an accident; it surrounds each secondary color by the primary
colors (pigments) out of which it can be made.  A puzzle composed by the first author asks whether it can be
recolored in such a way that {\em every} sector's color is a blend of the colors of that sector's neighbor,
without all the sectors being the same color.

It takes only a bit of tinkering to show that the answer is no, and (perhaps surprisingly) that no information about the
set of all colors or about how one might predict the color that arises from blending seems to be needed. Indeed, this
information might well depend on what brand of paint is in use. Instead, the only required assumption is
the following axiom about color.

\begin{axiom}[The Color Blend Axiom] \label{ax:blends}
The color that is created by a paint blend depends only on the amounts of each color in the blend (and not on the order in which
the colors were blended).
\end{axiom}

The color blend axiom is perhaps not {\em quite} as innocent as it looks, requiring, for example, that if we call an equal
mix of blue and yellow ``green," the blend must {\em behave} like green in other mixes.

If the puzzle were about numbers instead of colors, its solution (and generalization to edge-weighted, connected graphs)
would follow from the famous and extremely useful ``harmonic lemma" of graph theory, which says that if $f$ is a
real-valued function on the vertices of a finite, connected graph whose value at any vertex is the (weighted)
average of the values at its neighbors, then $f$ is constant. (See, e.g., Section 1.1.5 of \cite{D}.)

An obvious difficulty is that averaging is (arguably) the best-behaved function in the world, while paint blending
is a mystery (see Figure~\ref{fig:compare}).  In particular, the standard proof of the harmonic lemma
(by contradiction, from looking at a maximum-valued vertex whose neighbors are not all maximum-valued,
then invoking an order property of averaging) is unavailable for colors.

\begin{figure}
\centering
\renewcommand{\arraystretch}{1.3}
\begin{tabular}{|p{0.45\textwidth}|p{0.45\textwidth}|}
\hline
\textbf{Weighted Average} & \textbf{Paint Blending} \\
\hline
Maps continuously into an ordered field & Maps into a set \\
\hline
Has a nice formula & No formula known \\
\hline
``Linearity of expectation'' & Paint blending axiom \\
\hline
Respects order (arguments in an interval $\Rightarrow$ answer in the interval) & No order notion \\
\hline
\end{tabular}
\caption{An informal comparison.}\label{fig:compare}
\end{figure}

Nevertheless, we shall prove a ``chromonic lemma," in which the color wheel's 6-cycle is generalized to a finite,
edge-weighted, strongly connected digraph $G$.\footnote{$G$ is ``strongly connected" if for two vertices $u$ and
$v$ there is a directed path from $u$ to $v$.} We assume for convenience that $G$ has no loops, but in fact
loops can be easily accommodated as we will explain later in this section.

We can formalize the color blend axiom as follows.

Let $V$ be a finite set, and $\mathcal{B}$ the set of formal convex combinations (``blends") of the elements of $V$.
We regard $V$ as a subset of  $\mathcal{B}$ in the obvious way.

An equivalence relation $\equiv$ on $\mathcal{B}$ will be called a {\em color equivalence} if it respects
convex combination; that is, for every $x$, $y$, $z$ and $w \in \mathcal{B}$, and every $p$, $q \in [0,1]$,
$$
px + (1{-}p)y \equiv z \implies q(px + (1{-}p)y) + (1{-}q)w \equiv qz + (1{-}q)w.
$$
An equivalence class $C$ in $\mathcal{B}/\equiv$ is what we call a ``color"; $C(b)$ will denote the color of $b$,
i.e., the equivalence class to which the blend $b \in \mathcal{B}$ belongs.

Let $G = \langle V,E \rangle$ be a finite, edge-weighted, strongly connected digraph, and $\equiv$ a color
equivalence as above.  Thus, each $v \in V$ is a member of $\mathcal{B}$ as well as a vertex, and $C(v)$ is
that vertex's color.

If $(u,v)$ is an edge (from $u$ to $v$), its (positive, real) weight will be denoted by $w(u,v)$.
With $\equiv$ understood, $G$ will be called {\em chromonic} if for each $u \in V$,
$$
u \equiv \frac1{\sum_{i=1}^k w(u,v_i)} \left( \sum_{i=1}^k w(u,v_i)v_i \right)
$$
where $v_1, \dots, v_k$ are the out-neighbors of $u$.

\begin{figure}
\[
 \begin{tikzpicture}[>=stealth,scale=2, rotate =-90]
  \coordinate (a) at (-1,-1.73);
  \coordinate (b) at (1,-1.73);
  \coordinate (c) at (0,0);
  \coordinate (d) at (0,2);
  \draw (a) node[circle,draw=black,inner sep = 2] (A) {$u$};
  \draw (b) node[circle,draw=black,inner sep = 2] (B) {$v$};
  \draw (c) node[circle,draw=black,inner sep = 2] (C) {$w$};
  \draw (d) node[circle,draw=black,inner sep = 2] (D) {$x$};
  \draw[->] (A) to [bend right=15] node[midway,fill=white,draw=black,inner sep=1] {\footnotesize $1/2$} (B);
  \draw[->] (B) to [bend right=15] node[midway,fill=white,draw=black,inner sep=1] {\footnotesize $1/3$} (A);
  \draw[->] (A) to [bend right=15] node[midway,fill=white,draw=black,inner sep=1] {\footnotesize $1/2$} (C);
  \draw[->] (C) to [bend right=15] node[midway,fill=white,draw=black,inner sep=1] {\footnotesize $1/5$} (A);
  \draw[->] (B) to [bend right=15] node[midway,fill=white,draw=black,inner sep=1] {\footnotesize $2/3$} (C);
  \draw[->] (C) to [bend right=15] node[midway,fill=white,draw=black,inner sep=1] {\footnotesize $2/5$} (B);
  \draw[->] (C) to [bend right=15] node[midway,fill=white,draw=black,inner sep=1] {\footnotesize $2/5$} (D);
  \draw[->] (D) to [bend right=15] node[midway,fill=white,draw=black,inner sep=1] {\footnotesize $1$} (C);
 \end{tikzpicture}
\]
\caption{A strongly connected digraph $G$ with normalized edge weights.}\label{fig:digraph}
\end{figure}
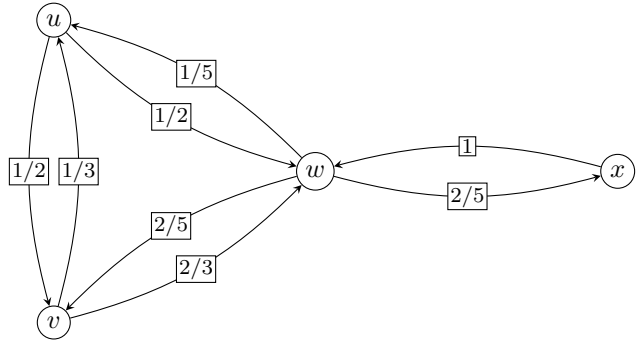
 
Note that there is no loss of generality if we normalize the edge-weights so that $\sum_{i=1}^k w(u,v_i) = 1$
whenever $u \in V$ with out-neighbors $v_1, \dots, v_k$, so we will henceforth assume $G$ is normalized,\footnote{To
normalize an undirected graph that doesn't happen to be regular, one must first replace each edge by two directed edges,
one in each direction.}  and the chromonicity condition then simplifies to
$$
u \equiv \sum_{i=1}^k w(u,v_i)v_i
$$
or, equivalently, to
$$
C(u) = C\left( \sum_{i=1}^k w(u,v_i)v_i \right)~.
$$

For the color wheel, $G$ is a $6$-cycle and the color of every vertex is a blend consisting of equal parts of its neighbors. 

In Figure \ref{fig:digraph} we see a more complicated example of a strongly connected digraph with normalized edge weights. The chromonicity condition says, for example, that color $u$ is a blend of equal parts color $v$ and color $w$, i.e., 
\[
C(u) = C\left( \frac{ v+w}{2}\right),
\] 
and color $w$ is a blend of one part $u$, two parts $v$, and two parts $x$, i.e.,
\[
C(w) = C\left( \frac{ u + 2v+2x}{5} \right).
\]
We also mention that from the point of view of blending, vertex $x$ must obviously have the same color as $w$ (since $w$ is its only neighbor). Thus the colors required for graph $G$ would be the same if we replace $x$ with a loop at $w$ of weight $2/5$. Conversely, this idea can be used to convert a chromonic digraph with loops to a loopless digraph without changing the set of colors.

Our ``chromonic lemma" is the following.

\begin{theorem}[The Chromonic Lemma]\label{thm:main}
If the finite, edge-weighted, strongly connected digraph $G$ is chromonic, all of its vertices have the same color.
\end{theorem}

\section{Proof Outline}

The general idea of the proof is to begin with a quantity of paint of color $C(v)$ at each vertex $v \in V$, and
distribute portions of this paint (broken down into the appropriate constituent colors) to neighbors of $v$
in such a way that the color of the paint at any vertex remains the same throughout. 

If for any $u$, $v \in V$ we can start with all the paint at $u$ and end with all the paint at $v$ after a
finite number of such steps, we will have shown that $C(u) = C(v)$ and proved the theorem.

Let $p_t(v)$ denote the amount of paint at vertex $v$ at time $t$, $t = 0,1,\dots$.  We may assume
$\sum_{v \in V} p_t(v) = 1$ for all $t$.  If $0 < h \le p_t(u)$, we may, at time $t$, execute a {\em push} of
magnitude $h$ at $u$ as follows: for each neighbor $v_i$ of $u$, an amount $hw(u,v_i)$ of paint of color $C(v_i)$
is sent from $u$ to $v_i$ (thus preserving the color at all vertices).  This is possible because $C(u) = \sum w(u,v_i)C(v_i)$.
If $h = p_t(u)$, so that $u$ is (temporarily) emptied of paint, the push is called a {\em shove}.

The inverse of a push, which we could call a ``pull," also has the desired properties and more closely corresponds to mixing real-world paint, but requires more conditions to ensure that the quantities $p_t(v_i)$ do not become negative.\footnote{It might help to imagine the push as normal paint mixing, but with time running in reverse.}

We can attempt to proceed by starting with $p_0(u) = 1$ and shoving vertices other than $v$ one by one until all
the paint is at $v$, but it takes infinite time to clear any pair of neighboring vertices, and we haven't got any
continuity assumptions on the space of colors to facilitate taking limits.

It suffices, however, to find a paint distribution $\pi$ that can be reached by finitely many pushes from the
all-paint-on-$u$ distribution for any $u$. To do this, it is both natural and convenient (but circumventable, at
some cost in complexity) to exploit a connection between paint distributions and random walks on $G$.

\section{Markov Chains and Stopping Rules}

Our normalized $G = \langle V,E \rangle$ corresponds to a unique irreducible finite-state Markov chain $\mathcal{M}$
with state space $V$ and transition matrix $P$ given by $P[u,v] = w(u,v)$. The chain $\mathcal{M}$ in turn has a unique stationary distribution $\pi$ satisfying $\pi P = P$.  (See, e.g., \cite{P}.)  This stationary distribution is a natural paint distribution to try to reach from any starting distribution, for two reasons:
(1) it has full support $V$, and (2) the problem of reaching $\pi$ (usually, approximately) by running the Markov
chain is important and very well studied indeed (see, e.g., \cite{H}).

The example in Figure \ref{fig:digraph} has
\[
 P = \left[ \begin{array}{rrrr}
  0 & \frac{1}{2} & \frac{1}{2} & 0 \\
  \frac{1}{3} & 0 & \frac{2}{3} & 0 \\
  \frac{1}{5} & \frac{2}{5} & 0 & \frac{2}{5} \\
  0 & 0 & 1 & 0
  \end{array}
 \right],
\]
and $\pi = (1/6, 1/4, 5/12, 1/6)$.

We think of $\mathcal{M}$ as the (edge-weighted) random walk of a token on $G$.  The amount $p_t(v)$ of paint at a vertex
$v$ is the probability that the token is at $v$ at time $t$. A shove at $v$ represents a step taken away from $v$, given
that the token is there; an $h$-push represents such a step taken with probability $h/p_t(v)$. A step by the token, wherever it might be, corresponds to simultaneous shoves at all the vertices of $G$ (we call this a {\em brawl}).  The number of steps it takes for the state distribution to reach $\pi$ {\em approximately} is called the ``mixing time" of $\mathcal{M}$. But we want to reach $\pi$ exactly, and for most Markov chains this never happens. However, since we don't need to restrict ourselves to brawls or even shoves, we are not limited to the token taking a fixed number of steps of the chain $\mathcal{M}$.

A {\em stopping rule} for $\mathcal{M}$ is an algorithm $A$ whose inputs are the current and past states of
the walk, and perhaps some random bits, to decide at any time whether to stop the walk or continue.

Four types of stopping rules, each stopping exactly at the stationary distribution (and employing the minimum
expected number of steps to do so), are described in \cite{LW1,LW2}.  Of these, one---the {\em threshold rule}---takes
a bounded number of steps.  It works as follows:  Each state $v$ is given a carefully chosen {\em threshold} $t_v$. (More on determining the threshold is in the next section.) If the walk reaches state $v$ at time $t \le t_v-1$, it continues; if at $t \ge t_v$, it terminates; if at $t \in (t_v-1, t_v)$, it continues with probability $t_v-t$.

We can easily duplicate the state distributions produced by the threshold rule, using paint distribution pushes.
We start at time 0 with a brawl, then at time 1 with another brawl, etc, but beginning at time $t_v$, vertex $v$ is no
longer pushed (but remains able to accept paint from neighbors). If $t_v$ is not an integer, then at time $t = \lfloor t_v \rfloor$,
vertex $v$ is not shoved but is pushed down only a fraction $t_v-t$ of its height.  By time $\max_v \lceil t_v \rceil$, we're done pushing and have reached the stationary paint distribution $\pi$.  Since this took only finitely many pushes, we have shown that the color obtained by blending the paints according to $\pi$ is $C(u)$; but $u$ was arbitrary, so all the colors are the same, and
the chromonic lemma is proved.

\section{Thresholds in action}

We don't need to compute the thresholds for our proof, but in fact the thresholds for a given $n$-state chain, starting distribution, and ending distribution, are computable in polynomial time. Here is one way to do it. Suppose we want to move from distribution $p=(p(v_1),\ldots,p(v_n))$ to distribution $q=(q(v_1),\ldots,q(v_n))$. Write $\Delta p = p-q$. (We will apply this to the case where $p$ is concentrated on a single vertex and $\pi$ is the stationary distribution.) Now identify a shove vector $h=(h(v_1),\ldots,h(v_n))$ as a solution to the linear system $h(P-I)=\Delta p$. The operator $P-I$ is known as the \emph{Laplacian} of $G$, and for strongly connected digraphs it has a one-dimensional nullspace spanned by $\pi$. That is, given a solution $h$, the vector $h+c\pi$ is also a solution for any scalar $c$. Thus we can choose $c$ so that the smallest entry of $h$ is zero. 

This now tells us exactly how much each vertex needs to be pushed to move from distribution $p$ to distribution $q$, and implicitly how many steps are required. We begin with $h_0(v) = h(v)$, and for $t=1,2,\ldots$ we do the following. If all $h_t(v)=0$ at time $t$, we stop. Otherwise, for each $v$, we push vertex $v$ by $\delta_t(v)/p_t(v)$, where $\delta_t(v)=\min\{ h_t(v), p_t(v) \}$, and set $h_{t+1}(v) = h_t(v)-\delta_t(v)$. For each $v$, there is a first time when $h_{t+1}(v) = 0$, and the threshold for $v$ is set to $t_v = t+\delta_t(v)/p_t(v)$.

\begin{figure}
\[
\begin{tikzpicture}[>=stealth]
\draw (0,0) node[scale=.75] {
\begin{tikzpicture}[>=stealth,scale=2, rotate =-90]
  \coordinate (a) at (-1,-1.73);
  \coordinate (b) at (1,-1.73);
  \coordinate (c) at (0,0);
  \coordinate (d) at (0,2);
  \draw (a) node[circle,draw=black,inner sep = 2] (A) {$0$};
  \draw (a) node[above left, xshift=-5] {\footnotesize \color{red}{$0$}};
  \draw (b) node[circle,draw=black,inner sep = 2] (B) {$0$};
  \draw (b) node[below left,xshift=-5] {\footnotesize \color{red}{$\frac{1}{20}$} };
  \draw (c) node[circle,draw=black,inner sep = 2] (C) {$0$};
  \draw (c) node[above,yshift=10] {\footnotesize \color{red}{$\frac{3}{4}$}};
  \draw (d) node[circle,draw=black,inner sep = 2] (D) {$1$};
  \draw (d) node[below right,xshift=5] {\footnotesize \color{red}{$\frac{17}{15}$}};
  \draw[->] (A) to [bend right=15] node[midway,fill=white,draw=black,inner sep=1] {\footnotesize $1/2$} (B);
  \draw[->] (B) to [bend right=15] node[midway,fill=white,draw=black,inner sep=1] {\footnotesize $1/3$} (A);
  \draw[->] (A) to [bend right=15] node[midway,fill=white,draw=black,inner sep=1] {\footnotesize $1/2$} (C);
  \draw[->] (C) to [bend right=15] node[midway,fill=white,draw=black,inner sep=1] {\footnotesize $1/5$} (A);
  \draw[->] (B) to [bend right=15] node[midway,fill=white,draw=black,inner sep=1] {\footnotesize $2/3$} (C);
  \draw[->] (C) to [bend right=15] node[midway,fill=white,draw=black,inner sep=1] {\footnotesize $2/5$} (B);
  \draw[->] (C) to [bend right=15] node[midway,fill=white,draw=black,inner sep=1] {\footnotesize $2/5$} (D);
  \draw[->] (D) to [bend right=15] node[midway,fill=white,draw=black,inner sep=1] {\footnotesize $1$} (C);
 \end{tikzpicture}
 };
 \draw (1,-4) node[scale=.75] {
\begin{tikzpicture}[>=stealth,scale=2, rotate =-90]
  \coordinate (a) at (-1,-1.73);
  \coordinate (b) at (1,-1.73);
  \coordinate (c) at (0,0);
  \coordinate (d) at (0,2);
  \draw (a) node[circle,draw=black,inner sep = 2] (A) {$0$};
  \draw (a) node[above left, xshift=-5] {\footnotesize \color{red}{$0$}};
  \draw (b) node[circle,draw=black,inner sep = 2] (B) {$0$};
  \draw (b) node[below left,xshift=-5] {\footnotesize \color{red}{$\frac{1}{20}$} };
  \draw (c) node[circle,draw=black,inner sep = 2] (C) {$1$};
  \draw (c) node[above,yshift=10] {\footnotesize \color{red}{$\frac{3}{4}$}};
  \draw (d) node[circle,draw=black,inner sep = 2] (D) {$0$};
  \draw (d) node[below right,xshift=5] {\footnotesize \color{red}{$\frac{2}{15}$}};
  \draw[->] (A) to [bend right=15] node[midway,fill=white,draw=black,inner sep=1] {\footnotesize $1/2$} (B);
  \draw[->] (B) to [bend right=15] node[midway,fill=white,draw=black,inner sep=1] {\footnotesize $1/3$} (A);
  \draw[->] (A) to [bend right=15] node[midway,fill=white,draw=black,inner sep=1] {\footnotesize $1/2$} (C);
  \draw[->] (C) to [bend right=15] node[midway,fill=white,draw=black,inner sep=1] {\footnotesize $1/5$} (A);
  \draw[->] (B) to [bend right=15] node[midway,fill=white,draw=black,inner sep=1] {\footnotesize $2/3$} (C);
  \draw[->] (C) to [bend right=15] node[midway,fill=white,draw=black,inner sep=1] {\footnotesize $2/5$} (B);
  \draw[->] (C) to [bend right=15] node[midway,fill=white,draw=black,inner sep=1] {\footnotesize $2/5$} (D);
  \draw[->] (D) to [bend right=15] node[midway,fill=white,draw=black,inner sep=1] {\footnotesize $1$} (C);
 \end{tikzpicture}
 };
 \draw (2,-8) node[scale=.75] {
\begin{tikzpicture}[>=stealth,scale=2, rotate =-90]
  \coordinate (a) at (-1,-1.73);
  \coordinate (b) at (1,-1.73);
  \coordinate (c) at (0,0);
  \coordinate (d) at (0,2);
  \draw (a) node[circle,draw=black,inner sep = 1] (A) {$\frac{3}{20}$};
  \draw (a) node[above left, xshift=-10] {\footnotesize \color{red}{$0$}};
  \draw (b) node[circle,draw=black,inner sep = 1] (B) {$\frac{3}{10}$};
  \draw (b) node[below left,xshift=-10] {\footnotesize \color{red}{$\frac{1}{20}$} };
  \draw (c) node[circle,draw=black,inner sep = 1] (C) {$\frac{1}{4}$};
  \draw (c) node[above,yshift=10] {\footnotesize \color{red}{$0$}};
  \draw (d) node[circle,draw=black,inner sep = 1] (D) {$\frac{3}{10}$};
  \draw (d) node[below right,xshift=10] {\footnotesize \color{red}{$\frac{2}{15}$}};
  \draw[->] (A) to [bend right=15] node[midway,fill=white,draw=black,inner sep=1] {\footnotesize $1/2$} (B);
  \draw[->] (B) to [bend right=15] node[midway,fill=white,draw=black,inner sep=1] {\footnotesize $1/3$} (A);
  \draw[->] (A) to [bend right=15] node[midway,fill=white,draw=black,inner sep=1] {\footnotesize $1/2$} (C);
  \draw[->] (C) to [bend right=15] node[midway,fill=white,draw=black,inner sep=1] {\footnotesize $1/5$} (A);
  \draw[->] (B) to [bend right=15] node[midway,fill=white,draw=black,inner sep=1] {\footnotesize $2/3$} (C);
  \draw[->] (C) to [bend right=15] node[midway,fill=white,draw=black,inner sep=1] {\footnotesize $2/5$} (B);
  \draw[->] (C) to [bend right=15] node[midway,fill=white,draw=black,inner sep=1] {\footnotesize $2/5$} (D);
  \draw[->] (D) to [bend right=15] node[midway,fill=white,draw=black,inner sep=1] {\footnotesize $1$} (C);
 \end{tikzpicture}
 };
\draw (3,-12) node[scale=.75] {
\begin{tikzpicture}[>=stealth,scale=2, rotate =-90]
  \coordinate (a) at (-1,-1.73);
  \coordinate (b) at (1,-1.73);
  \coordinate (c) at (0,0);
  \coordinate (d) at (0,2);
  \draw (a) node[circle,draw=black,inner sep = 0] (A) {$\frac{1}{6}$};
  \draw (a) node[above left, xshift=-10] {\footnotesize \color{red}{$0$}};
  \draw (b) node[circle,draw=black,inner sep = 0] (B) {$\frac{1}{4}$};
  \draw (b) node[below left,xshift=-10] {\footnotesize \color{red}{$0$} };
  \draw (c) node[circle,draw=black,inner sep = 0] (C) {$\frac{5}{12}$};
  \draw (c) node[above,yshift=10] {\footnotesize \color{red}{$0$}};
  \draw (d) node[circle,draw=black,inner sep = 0] (D) {$\frac{1}{6}$};
  \draw (d) node[below right,xshift=10] {\footnotesize \color{red}{$0$}};
  \draw[->] (A) to [bend right=15] node[midway,fill=white,draw=black,inner sep=1] {\footnotesize $1/2$} (B);
  \draw[->] (B) to [bend right=15] node[midway,fill=white,draw=black,inner sep=1] {\footnotesize $1/3$} (A);
  \draw[->] (A) to [bend right=15] node[midway,fill=white,draw=black,inner sep=1] {\footnotesize $1/2$} (C);
  \draw[->] (C) to [bend right=15] node[midway,fill=white,draw=black,inner sep=1] {\footnotesize $1/5$} (A);
  \draw[->] (B) to [bend right=15] node[midway,fill=white,draw=black,inner sep=1] {\footnotesize $2/3$} (C);
  \draw[->] (C) to [bend right=15] node[midway,fill=white,draw=black,inner sep=1] {\footnotesize $2/5$} (B);
  \draw[->] (C) to [bend right=15] node[midway,fill=white,draw=black,inner sep=1] {\footnotesize $2/5$} (D);
  \draw[->] (D) to [bend right=15] node[midway,fill=white,draw=black,inner sep=1] {\footnotesize $1$} (C);
 \end{tikzpicture}
 };
 \draw[->] (1,-1) to [bend left=15] node[midway, right] {Shove $x$} (1.5,-3);
 \draw[->] (2,-5) to [bend left=15] node[midway, right] {$\frac{3}{4}$-push of $w$} (2.5,-7);
 \draw[->] (3,-9) to [bend left=15] node[midway, right] {$\frac{1}{6}$-push of $v$, $\frac{4}{9}$-push of $x$,} (3.5,-11);
 \end{tikzpicture}
\]
\caption{Pushing on the graph to move paint from $x$ to the stationary distribution. Total push values are listed next to each vertex in red.}\label{fig:pushes}
\end{figure}
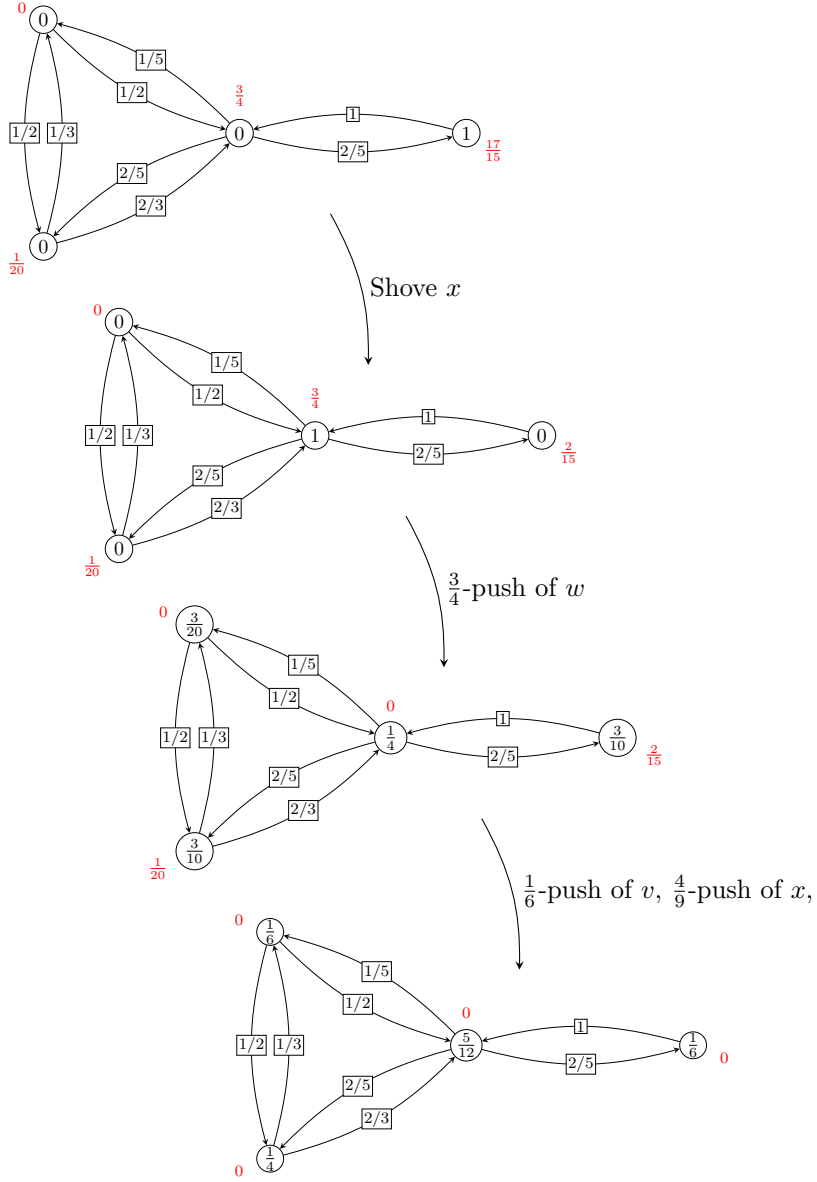

For the running example, if we wish to begin with all the paint on vertex $x$, this process leads us to find the following total push vector: $h=(0,1/20, 3/4, 17/15)$. We can apply one shove and three pushes to move from $p_0(x)=1$ to the stationary distribution $\pi$. The thresholds are: $(0, 13/6, 7/4, 22/9)$. See Figure \ref{fig:pushes}.

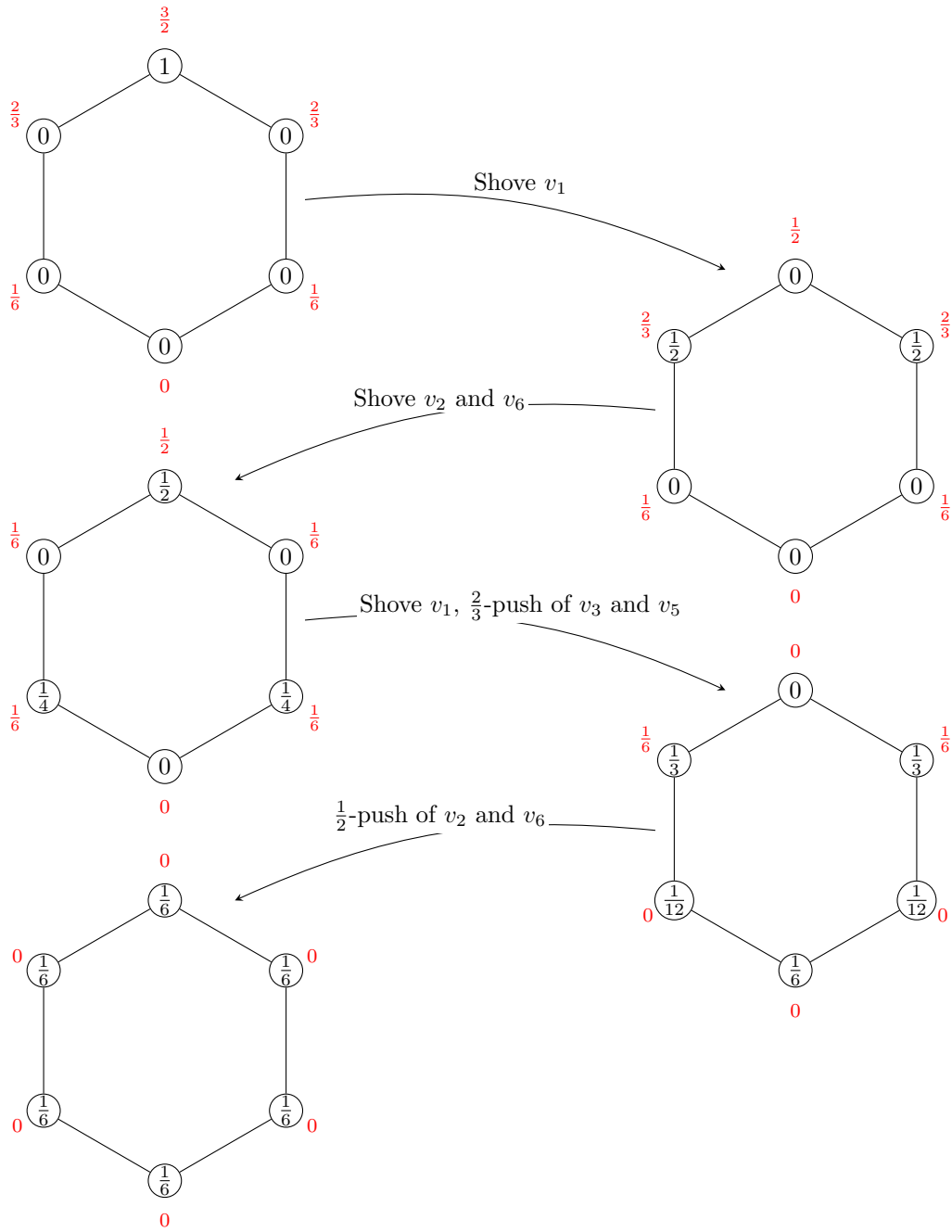
\begin{figure}
\[
\begin{tikzpicture}[>=stealth]
\draw (0,0) node {
\begin{tikzpicture}[scale=1]
  \coordinate (a) at (0,2);
  \coordinate (b) at (-1.73,1);
  \coordinate (c) at (-1.73,-1);
  \coordinate (d) at (0,-2);
  \coordinate (e) at (1.73,-1);
  \coordinate (f) at (1.73,1);
  \draw (a) node[circle,draw=black,inner sep = 2] (A) {$1$};
  \draw (a) node[above, yshift=10] {\footnotesize \color{red}{$\frac{3}{2}$}};
  \draw (b) node[circle,draw=black,inner sep = 2] (B) {$0$};
  \draw (b) node[above left,xshift=-5] {\footnotesize \color{red}{$\frac{2}{3}$} };
  \draw (c) node[circle,draw=black,inner sep = 2] (C) {$0$};
  \draw (c) node[below left,xshift=-5] {\footnotesize \color{red}{$\frac{1}{6}$}};
  \draw (d) node[circle,draw=black,inner sep = 2] (D) {$0$};
  \draw (d) node[below,yshift=-10] {\footnotesize \color{red}{$0$}};
  \draw (e) node[circle,draw=black,inner sep = 2] (E) {$0$};
  \draw (e) node[below right,xshift=5] {\footnotesize \color{red}{$\frac{1}{6}$}};
  \draw (f) node[circle,draw=black,inner sep = 2] (F) {$0$};
  \draw (f) node[above right,xshift=5] {\footnotesize \color{red}{$\frac{2}{3}$}};
  \draw (A) -- (B) -- (C)--(D)--(E)--(F)--(A);
 \end{tikzpicture}
 };
 \draw (9,-3) node {
\begin{tikzpicture}[scale=1]
  \coordinate (a) at (0,2);
  \coordinate (b) at (-1.73,1);
  \coordinate (c) at (-1.73,-1);
  \coordinate (d) at (0,-2);
  \coordinate (e) at (1.73,-1);
  \coordinate (f) at (1.73,1);
  \draw (a) node[circle,draw=black,inner sep = 2] (A) {$0$};
  \draw (a) node[above, yshift=10] {\footnotesize \color{red}{$\frac{1}{2}$}};
  \draw (b) node[circle,draw=black,inner sep = 0] (B) {$\frac{1}{2}$};
  \draw (b) node[above left,xshift=-5] {\footnotesize \color{red}{$\frac{2}{3}$} };
  \draw (c) node[circle,draw=black,inner sep = 2] (C) {$0$};
  \draw (c) node[below left,xshift=-5] {\footnotesize \color{red}{$\frac{1}{6}$}};
  \draw (d) node[circle,draw=black,inner sep = 2] (D) {$0$};
  \draw (d) node[below,yshift=-10] {\footnotesize \color{red}{$0$}};
  \draw (e) node[circle,draw=black,inner sep = 2] (E) {$0$};
  \draw (e) node[below right,xshift=5] {\footnotesize \color{red}{$\frac{1}{6}$}};
  \draw (f) node[circle,draw=black,inner sep = 0] (F) {$\frac{1}{2}$};
  \draw (f) node[above right,xshift=5] {\footnotesize \color{red}{$\frac{2}{3}$}};
  \draw (A) -- (B) -- (C)--(D)--(E)--(F)--(A);
 \end{tikzpicture}
 };
 \draw (0,-6) node {
\begin{tikzpicture}[scale=1]
  \coordinate (a) at (0,2);
  \coordinate (b) at (-1.73,1);
  \coordinate (c) at (-1.73,-1);
  \coordinate (d) at (0,-2);
  \coordinate (e) at (1.73,-1);
  \coordinate (f) at (1.73,1);
  \draw (a) node[circle,draw=black,inner sep = 0] (A) {$\frac{1}{2}$};
  \draw (a) node[above, yshift=10] {\footnotesize \color{red}{$\frac{1}{2}$}};
  \draw (b) node[circle,draw=black,inner sep = 2] (B) {$0$};
  \draw (b) node[above left,xshift=-5] {\footnotesize \color{red}{$\frac{1}{6}$} };
  \draw (c) node[circle,draw=black,inner sep = 0] (C) {$\frac{1}{4}$};
  \draw (c) node[below left,xshift=-5] {\footnotesize \color{red}{$\frac{1}{6}$}};
  \draw (d) node[circle,draw=black,inner sep = 2] (D) {$0$};
  \draw (d) node[below,yshift=-10] {\footnotesize \color{red}{$0$}};
  \draw (e) node[circle,draw=black,inner sep = 0] (E) {$\frac{1}{4}$};
  \draw (e) node[below right,xshift=5] {\footnotesize \color{red}{$\frac{1}{6}$}};
  \draw (f) node[circle,draw=black,inner sep = 2] (F) {$0$};
  \draw (f) node[above right,xshift=5] {\footnotesize \color{red}{$\frac{1}{6}$}};
  \draw (A) -- (B) -- (C)--(D)--(E)--(F)--(A);
 \end{tikzpicture}
 };
 \draw (9,-9) node {
\begin{tikzpicture}[scale=1]
  \coordinate (a) at (0,2);
  \coordinate (b) at (-1.73,1);
  \coordinate (c) at (-1.73,-1);
  \coordinate (d) at (0,-2);
  \coordinate (e) at (1.73,-1);
  \coordinate (f) at (1.73,1);
  \draw (a) node[circle,draw=black,inner sep = 2] (A) {$0$};
  \draw (a) node[above, yshift=10] {\footnotesize \color{red}{$0$}};
  \draw (b) node[circle,draw=black,inner sep = 0] (B) {$\frac{1}{3}$};
  \draw (b) node[above left,xshift=-5] {\footnotesize \color{red}{$\frac{1}{6}$} };
  \draw (c) node[circle,draw=black,inner sep = 0] (C) {$\frac{1}{12}$};
  \draw (c) node[below left,xshift=-5] {\footnotesize \color{red}{$0$}};
  \draw (d) node[circle,draw=black,inner sep = 0] (D) {$\frac{1}{6}$};
  \draw (d) node[below,yshift=-10] {\footnotesize \color{red}{$0$}};
  \draw (e) node[circle,draw=black,inner sep = 0] (E) {$\frac{1}{12}$};
  \draw (e) node[below right,xshift=5] {\footnotesize \color{red}{$0$}};
  \draw (f) node[circle,draw=black,inner sep = 0] (F) {$\frac{1}{3}$};
  \draw (f) node[above right,xshift=5] {\footnotesize \color{red}{$\frac{1}{6}$}};
  \draw (A) -- (B) -- (C)--(D)--(E)--(F)--(A);
 \end{tikzpicture}
 };
 \draw (0,-12) node {
\begin{tikzpicture}[scale=1]
  \coordinate (a) at (0,2);
  \coordinate (b) at (-1.73,1);
  \coordinate (c) at (-1.73,-1);
  \coordinate (d) at (0,-2);
  \coordinate (e) at (1.73,-1);
  \coordinate (f) at (1.73,1);
  \draw (a) node[circle,draw=black,inner sep = 0] (A) {$\frac{1}{6}$};
  \draw (a) node[above, yshift=10] {\footnotesize \color{red}{$0$}};
  \draw (b) node[circle,draw=black,inner sep = 0] (B) {$\frac{1}{6}$};
  \draw (b) node[above left,xshift=-5] {\footnotesize \color{red}{$0$} };
  \draw (c) node[circle,draw=black,inner sep = 0] (C) {$\frac{1}{6}$};
  \draw (c) node[below left,xshift=-5] {\footnotesize \color{red}{$0$}};
  \draw (d) node[circle,draw=black,inner sep = 0] (D) {$\frac{1}{6}$};
  \draw (d) node[below,yshift=-10] {\footnotesize \color{red}{$0$}};
  \draw (e) node[circle,draw=black,inner sep = 0] (E) {$\frac{1}{6}$};
  \draw (e) node[below right,xshift=5] {\footnotesize \color{red}{$0$}};
  \draw (f) node[circle,draw=black,inner sep = 0] (F) {$\frac{1}{6}$};
  \draw (f) node[above right,xshift=5] {\footnotesize \color{red}{$0$}};
  \draw (A) -- (B) -- (C)--(D)--(E)--(F)--(A);
 \end{tikzpicture}
 };
 \draw[->] (2,0) to [bend left=15] node[midway, above] {Shove $v_1$} (8,-1);
 \draw[->] (7,-3) to [bend right=15] node[midway, above,fill=white,inner sep=1] {Shove $v_2$ and $v_6$} (1,-4);
 \draw[->] (2,-6) to [bend left=15] node[midway, above, fill=white,inner sep=1] {Shove $v_1$, $\frac{2}{3}$-push of $v_3$ and $v_5$} (8,-7);
 \draw[->] (7,-9) to [bend right=15] node[midway, above, fill=white,inner sep=1] {$\frac{1}{2}$-push of $v_2$ and $v_6$} (1,-10);
 \end{tikzpicture}
\]
\caption{The threshold rule in action on the graph $C_6$.}\label{fig:push2}
\end{figure}

For a simpler example, we can return to the color wheel itself, where $G$ is the undirected 6-cycle with equal
edge-weights.  For this graph, the stationary distribution is uniform (as it must be, by symmetry). Supposing we start with all the paint on vertex $1$, we find the total push vector to be $h=(\frac{3}{2},\frac{2}{3}, \frac{1}{6}, 0, \frac{1}{6}, \frac{2}{3})$.  The thresholds are $(3, \frac{7}{2}, \frac{8}{3}, 0, \frac{8}{3}, \frac{7}{2})$. In Figure \ref{fig:push2} we see the shoves and pushes that take the paint from one vertex to the stationary (uniform) distribution. We conclude that the paint initially at vertex 1 is the same color as the equal blend of the colors at all vertices,
and by symmetry, the same as the color at every other vertex.

\end{document}